\documentstyle[12pt]{article}

\newtheorem{theorem}{Theorem}

\newtheorem{proposition}{Proposition}
\newtheorem{definition}{Definition}

\newtheorem{remark}{Remark}

\begin{document}

\title{Twisting Elements in  Homotopy G-algebras}

\author{T. Kadeishvili}

\date{ }

\maketitle

\begin{abstract}
 We study the notion of twisting elements $da=a\smile_1a$ with
 respect to $\smile_1$ product when it is a part of homotopy
 Gerstenhaber algebra structure.
 This allows to bring to one context the two classical concepts,
 the theory of deformation of algebras of M. Gerstenhaber, and
 $A(\infty)$-algebras of J. Stasheff.
 \end{abstract}


\section{Introduction}
\label{intro}

A \emph{twisting element} in a differential graded algebra (dga)
$(A=\{A^i\},d:A^n\to A^{n+1},a^m\cdot b^n\in A^{m+n})$ is defined
as an element $t\in A^1$ satisfying the Brown's condition
\begin{equation}
\label{brow} dt=t\cdot t.
\end{equation}
Denote the set of all twisting elements by $Tw(A)$. An useful
consequences of the Brown' condition is the following: let $M$ be
a dg module over $A$, then a twisting element  $t \in Tw(A)$
defines on $M$ a new differential $d_t:M\to M$  by
$d_t(x)=dx+t\cdot x$, and the condition (\ref{brow}) guarantees
that $d_td_t=0$.

Twisting elements show up in various problems of algebraic
topology and homological algebra. The first appearance was in
homology theory of fibre bundles \cite{Bro59}: For a fibre bundle
$F\to E\to B$ with structure group $G$ there exists a twisting
element $t\in A=C^*(B,C_*(G))$ such that $(M=C_*(B)\otimes
C_*(F),d_t)$ (the twisted tensor product) gives homology of the
total space $E$.

Later N. Berikashvili \cite{Ber76} has introduced in $Tw(A)$ an
\emph{equivalence relation} induced by the following  group
action. Let $G$ be the group of invertible elements in $A^0$, then
for $g\in G$ and $t\in Tw(A)$ let
\begin{equation}
\label{berik} g\ast t=g\cdot t\cdot g^{-1}+dg\cdot g^{-1},
\end{equation}
easy to see that $g\ast t\in Tw(A)$. The factor set
$D(A)=Tw(A)/G$, called Berikashvili's functor $D$, has nice
properties and useful applications. In particular if $t\sim t'$
then $(M,d_t)$ and $(M,d_{t'})$ are {\it isomorphic}.

The notion of \emph{homotopy G-algebra} (hGa in short) was
introduced by Gerstenhaber and Voronov in \cite{GV95} as an
additional structure on a dg algebra $(A,d,\cdot)$ that induces a
Gerstenhaber algebra structure on homology. The main example is
Hochschild cochain complex of an algebra.

Another point of view is that hGa is a particular case of
$B(\infty)$-algebra \cite{GJ94}: this is an additional structure
on a dg algebra $(A,d,\cdot)$ that induces a dg bialgebra
structure on the bar construction $BA$.

There is the third aspect of hGa \cite{Kad04}: this is a structure
which measures the noncommutativity of $A$. The Steenrod's
$\smile_1$ product which is the classical tool which measures the
noncommutativity of a dg algebra $(A,d,\cdot )$  satisfies the
condition
\begin{equation}
\label{Steen} d(a\smile_1 b)=da\smile_1 b + a\smile_1 db +a\cdot b
+ b\cdot a.
\end{equation}
The existence of such $\smile_1$ guarantees the commutativity of
$H(A)$, but a $\smile _1$ product satisfying just the condition
(\ref{Steen}) is too poor for some applications. In many
constructions some deeper properties of $\smile_1$ are needed, for
example the compatibility with the dot product of $A$ (the Hirsch
formula)
\begin{equation}
\label{Hirsch} (a\cdot b)\smile _1c + a\cdot (b\smile
_1c)+(a\smile _1c)\cdot b=0.
\end{equation}
A hGa $(A,d,\cdot,\{E_{1,k}\})$ is a dga $(A,d,\cdot)$ equipped
additionally with a sequence of operations (some authors call them
\emph{braces})
$$
\{E_{1,k}:A\otimes A^{\otimes k}\to A,\ k=1,2,...\}
$$
satisfying some coherence conditions (see bellow). The starting
operation $E_{1,1}$ is a kind of $\smile_1$ product: it satisfies
the conditions (\ref{Steen}) and (\ref{Hirsch}). As for the
symmetric expression
$$
a\smile _1(b\cdot c)+ b\cdot (a\smile _1c)+(a\smile _1b)\cdot c,
$$
it is just {\it homotopical to zero} and the appropriate chain
homotopy is the operation $E_{1,2}$. So we can say that a hGa is a
dga with a "good" $\smile_1$ product.

There is one more aspect of hGa: the operation $E_{1,1}=\smile_1$
is not associative but the commutator
$[a,b]=a\smile_1b-b\smile_1a$ satisfies the Jacobi identity, so it
forms on the desuspension $s^{-1}A$ a structure of dg Lie algebra.

Let us present three remarkable examples of homotopy G-algebras.

The first one is the cochain complex of 1-reduced simplicial set
$C^*(X)$. The operations $E_{1,k}$ here are dual to cooperations
defined by Baues in \cite{Bau81}, and the starting operation
$E_{1,1}$ is  the classical Steenrod's $\smile_1$ product.

The second example is the Hochschild cochain complex $C^*(U,U)$ of
an associative algebra $U$. The operations $E_{1,k}$ here were
defined in \cite{Kad88} with the purpose to describe
$A(\infty)$-algebras in terms of Hochschild cochains although the
properties of those operations which where used as defining ones
for the notion of homotopy G-algebra in \cite{GV95} did not appear
there. These operations where defined also in \cite{GJ94}. Again
the starting operation $E_{1,1}$ is  the classical Gerstenhaber's
circle product which is sort of $\smile_1$-product.

The third example is  the the cobar construction $\Omega C$ of a
dg \emph{bialgebra} $C$. The operations $E_{1,k}$ are constructed
in \cite{Kad05}. And again the starting operation $E_{1,1}$ is
classical: it is  Adams's $\smile_1$-product defined for $\Omega
C$ in \cite{Ada60} using the \emph{multiplication} of $C$.

The main task of this paper is to introduce the notion of a
twisting  element and their transformation in a hGa. Shortly a
twisting element in a hGa $(A,d,\cdot,\{E_{1,k}\})$ is an element
$a\in A$ such that $da=a\smile_1a$ and two twisting elements
$a,\overline{a}\in A$ we call equivalent if there exists $g\in A$
such that
$$
 \overline{a}=a+dg+g\cdot
g+g\smile_1a+\overline{a}\smile_1g+E_{1,2}(\overline{a};g,g)+
E_{1,3}(\overline{a};g,g,g)+...\ .
$$
As we see in the definition of a twisting element participates
just the operation $E_{1,1}=\smile_1$ but in the definition of
equivalence participates the whole hGa structure. We remark that
such a twisting element $a\in A$ is a Lie twisting element in the
dg Lie algebra $(s^{-1}A,d,[\ ,\ ])$, i.e. satisfies
$da=\frac{1}{2}[a,a]$. But it is unclear wether the equivalence
can be formulated in terms the bracket $[\ ,\ ]$.

In this paper we present the following application of the notion
of twisting element in a hGa: it allows to unify two classical
concepts, namely the theory of deformation of algebras of M.
Gerstenhaber, and  J. Stasheff's $A(\infty)$-algebras.

Namely, a Gerstenhaber's deformation of an associative algebra $U$
(see \cite{Ger64}, and bellow)
$$
a\star b=a\cdot b+B_1(a\otimes b)t+B_2(a\otimes b)t^2+B_3(a\otimes
b)t^3+...\in U[[t]],
$$
can be considered as a twisting element $B=B_1+B_2+...\in
C^2(U,U)$ in the Hochschild cochain complex of $U$ with
coefficients in itself: the defining condition of deformation
means exactly $dB=B\smile_1B$. Furthermore, two deformations are
equivalent if and only if the corresponding twisting elements are
equivalent in the above sense.

On the other side, the same concept of twisting elements in hGa
works in the following problem. Suppose $(H,\mu :H\otimes
H\rightarrow H)$ is a graded algebra. Let us define it's \emph{
Stasheff's deformation} as an $A(\infty)$ algebra structure
$(H,\{m_i\})$ with $m_1=0$ and $m_2=\mu$, i.e.  which extends the
given algebra structure. Then each deformation can be considered
as a twisting element $m=m_3+m_4+...,\ m_i\in C^i(H,H)$ in the
Hochschild cochain complex of $H$ with coefficients in itself: the
Stasheffs defining condition of $A(\infty)$-algebra means exactly
$dm=m\smile_1m$. Furthermore, to isomorphic  (as
$A(\infty)$-algebras) deformations correspond equivalent twisting
elements in the above sense.

In both cases we present the obstruction theory for the existence
of suitable deformations. The obstructions live in suitable
Hochschild cohomologies: in $H^2(U,U)$ in Gerstenhaber's
deformation case and in $H^i(H,H),\ i=3,4,...$ in Stasheff's
deformation case.

The structure of the paper is following. In the section
\ref{notations} necessary definitions are given. In the section
\ref{HGA} the definition Homotopy G-algebra is presented. In the
section \ref{twisting} the notion of twisted element in a homotopy
G-algebra is studied. In the last two sections \ref{gerdef} and
\ref{stadef} the above mentioned applications of this notion are
given.

\noindent {\bf Acknowledgements.} Dedicated to Murray
Gerstenhaber's 80th and Jim Stasheff's 70th birthdays.

\section{Notation and Preliminaries}
\label{notations}

We work over $R=Z_2$. For a graded module $M$ we denote by $sM$
the suspension of $M$, i.e. $(sM)^i=M^{i-1}$ . Respectively
$(s^{-1}M)^i=M^{i+1}$.

\subsection{Differential Graded Algebras and Coalgebras}
A {\em differential graded algebra } (dg algebra, or dga) is a
graded R-module $C=\{C^{i},\ i\in
 Z\}$ with an associative \emph{multiplication} $\mu :C^i\otimes
C^j\to  C^{i+j}$ and a \emph{differential} $d:C^i \to C^{i+1}$
satisfying $dd=0$ and the \emph{derivation condition} $d(x\cdot
y)=dx\cdot y+x\cdot dy$, where $x\cdot y=\mu(x\otimes y)$. A dga
$C$ is {\em connected} if $C^{<0}=0$ and $C^0=R.$ A connected dga
$C$ is {\em n-reduced} if $C^i=0$ for $1\leq i\leq n$.

 A {\em differential graded
coalgebra } (dg coalgebra, or dgc) is a graded $R$-module
$C=\{C^{i},\ i\in Z\}$ with a coassociative
\emph{comultiplication} $\Delta :C\to C\otimes C$ and a
\emph{differential} $d:C^i \to C^{i+1}$ satisfying $dd=0$ and the
\emph{coderivation condition}$\Delta d=(d\otimes id +id\otimes
d)\Delta$. A dgc $C$ is {\em connected} if $C_{<0}=0$ and $C_0=R$.
A connected  dgc $C$ is {\em n-reduced} if $C_i=0$ for $1\leq
i\leq n$.

A {\em differential graded  bialgebra } (dg bialgebra) $(C,d,\mu ,
\Delta)$ is a dg coalgebra $(C,d,\Delta)$ with a morphism of dg
coalgebras $\mu :C\otimes C\to  C$ turning $(C,d,\mu )$ into a dg
algebra.

\subsection{Cobar and Bar Constructions}
\label{bar}

Let $M$ be a graded $R$-module with $M^{i\leq 0}=0$ and let $T(M)$
be the tensor algebra of $M$, i.e. $ T(M) =\oplus_{i=0}^{\infty}
M^{\otimes i}$. Tensor algebra $T(M)$ is a free graded algebra
generated by $M$: for a graded algebra $A$ and a homomorphism
$\alpha:M\to A$ of degree zero there exists a unique morphism of
graded algebras $f_{\alpha}:T(M)\to A$ (called {\it multiplicative
extension} of $\alpha$)such that $f_{\alpha}(a)=\alpha(a)$. The
map $f_{\alpha}$ is given by $ f_{\alpha}(a_1\otimes...\otimes
a_n)=\alpha(a_1)\cdot ... \cdot \alpha(a_n)$.

Dually, let $T^c(M)$ be the tensor coalgebra of $M$, i.e. $ T^c(M)
=\oplus_{i=0}^{\infty}  M^{\otimes i}$, and the comultiplication
$\nabla :  T^c(M)\to  T^c(M)\otimes  T^c(M)$ is given by
$$
\nabla (a_1\otimes...\otimes a_n)=\sum_{k=0}^n
(a_1\otimes...\otimes a_k)\otimes (a_{k+1}\otimes...\otimes a_n).
$$
The tensor coalgebra $(T^c(M),\nabla)$ is a cofree graded
coalgebra: for a graded coalgebra $C$ and a homomorphism
$\beta:C\to M$ of degree zero there exists a unique morphism of
graded coalgebras $g_{\beta}:C\to T^c(M)$ (called {\it
comultiplicative coextension} of $\beta$) such that
$p_1g_{\beta}=\beta $, here $p_1:T^c(M)\to M$ is the clear
projection. The map $g_{\beta}$ is given by
$$
g_{\beta}=\sum_{n=0}^{\infty}
(\beta\otimes...\otimes\beta)\Delta^n,
$$
where $\Delta^n:C\to C^{\otimes n}$ is  $n$-th iteration of the
diagonal $\Delta:C\to C\otimes C$, i.e. $\Delta^1=id,\
\Delta^2=\Delta,\ \Delta^n=(\Delta^{n-1} \otimes id)\Delta$.

Let $(C,d_C, \Delta )$ be a connected dg coalgebra and $\Delta(c)
=c\otimes 1_R+1_R\otimes c+   \Delta' (c)$. The (reduced) {\it
cobar construction} $ \Omega C$ on $C$ is a dg algebra whose
underlying graded algebra is $T(sC^{>0})$. An element
$(sc_1\otimes...\otimes sc_n)\in (sC)^{\otimes  n}\subset
T(sC^{>0})$ is denoted by $[c_1,...,c_n]\in \Omega C$. The
differential $d_\Omega$ of $\Omega C$ for a generator $[c]\in
\Omega C $ is defined by $d_\Omega[c]=[d_C(c)]+ \sum [c',c'']$
where $\Delta'(c)=\sum c'\otimes c''$, and is extended as a
derivation.

Let $(A,d_A, \mu)$ be a 1-reduced  dg algebra. The (reduced) {\it
bar construction} $ BA$ on $A$ is a dg coalgebra whose underlying
graded coalgebra is $T^c(s^{-1}A^{>0})$. Again an element
$(s^{-1}a_1\otimes...\otimes s^{-1}a_n)\in (s^{-1}A)^{\otimes
n}\subset T^c(s^{-1}A^{>0})$ we denote as $[a_1,...,a_n]\in BA$.
The differential $d_B$ of $BA$ is defined by
$$
\begin{array}{l}
d_B[a_1,...,a_n]=\sum_{i=1}^{n} [a_1,...,d_Aa_i,...,a_n]
+\sum_{i=1}^{n-1} [a_1,...,a_{i}\cdot a_{i+1},...,a_n].
\end{array}
$$


\subsection{Twisting Cochains}
\label{twistcoch}

Let $(C,d,\Delta)$ be a dgc and $(A,d,\mu)$ be a dga. A twisting
cochain \cite{Bro59} is a homomorphism $\tau:C\to A$ of degree +1
satisfying the Brown's condition
\begin{equation}
\label{Bro59} d\tau+\tau d=\tau\smile \tau,
\end{equation}
where $\tau\smile \tau'=\mu_A(\tau\otimes\tau')\Delta $. We denote
by $Tw(C,A)$ the set of all twisting cochains $\tau:C\to A$.

There are universal twisting cochains $\tau_C:C\to \Omega C$ and
$\tau_A:BA\to A$ being clear inclusion and projection
respectively.

Here are essential consequences of the condition (\ref{Bro59}):

\noindent (i) The multiplicative extension $f_{\tau}:\Omega C\to
A$ is a map of dg algebras, so there is a bijection
$Tw(C,A)\leftrightarrow Hom_{dg-Alg}(\Omega C,A)$;

\noindent (ii) The comultiplicative coextension $g_{\tau}:C\to BA$
is a map of dg coalgebras, so there is a bijection
$Tw(C,A)\leftrightarrow Hom_{dg-Coalg}(C,BA)$.


\section{Homotopy G-algebras}
\label{HGA}

A {\it homotopy G-algebra } (hGa in short) is a dg algebra with
"good" $\smile_1$ product. The general notion was introduced in
\cite{GV95}, see also \cite{Vor99}.

\begin{definition} A homotopy G-algebra is defined as a dg algebra $(A,d,\cdot )$ with a
given sequence of operations
$$
E_{1,k}:A\otimes (A^{\otimes k})\rightarrow A,\quad k=0,1,2,3,...
$$
(the value of the operation $E_{1,k}$ on $a\otimes b_1\otimes
...\otimes b_k\in A\otimes (A\otimes ...\otimes A)$ we write as
$E_{1,k}(a;b_1,...,b_k)$) which satisfies the conditions
\begin{equation}
\label{E10}
\begin{array}{l}
E_{1,0}=id,
\end{array}
\end{equation}

\begin{equation}
\label{E1n}
\begin{array}{l}
dE_{1,k}(a;b_1, ..., b_k)+ E_{1,k}(da;b_1, ...,b_k)+ \sum_i
E_{1,k}(a;b_1, ..., db_i, ..., b_k)\\
=  b_1\cdot E_{1,k-1}(a;b_2, ..., b_k) +
E_{1,k-1}(a;b_1, ..., b_{k-1})\cdot b_k+\\
\ \ \ \sum_iE_{1,k-1}(a;b_1, ..., b_i\cdot b_{i+1}, ..., b_k),
\end{array}
\end{equation}
\begin{equation}
\label{E2n}
\begin{array}{l}
E_{1,k}(a_1\cdot a_2;b_1,.., b_k)\\
=
a_1\cdot E_{1,k}(a_2;b_1, ..., b_k)+ E_{1,k}(a_1;b_1, ..., b_k)\cdot a_2+\\
\ \ \ \sum_{p=1}^{k-1} E_{1,p}(a_1;b_1, ..., b_p)\cdot
E_{1,m-p}(a_2;b_{p+1}, ..., b_k),
\end{array}
\end{equation}
\begin{equation}
\label{E1assoc}
\begin{array}{l}
E_{1,n}(E_{1,m}(a ;b _1, ..., b _m);c_1, ..., c _n)\\
= \sum_{0\leq
i_1\leq j_1\leq...\leq i_m\leq j_m\leq n}
\\
\ \ \ E_{1,n-(j_1+...+j_m)+(i_1+...+i_m)+m}(a ;c _1, ..., c
_{i_1},E_{1,j_1-i_1}(b_1;c _{i_1+1}, ..., c _{j_1}), \\
\ \ \ c _{j_1+1}, ...,  c _{i_2},E_{1,j_2-i_2}(b _2;c _{i_2+1},
..., c _{j_2}), c _{j_2+1}, ..., c _{i_m},\\
\ \ \ E_{1,j_m-i_m}(b _m;c _{i_m+1}, ..., c _{j_m}), c
_{j_m+1},...,c _n).
\end{array}
\end{equation}
\end{definition}


Let us present these conditions in low dimensions.

The condition (\ref{E1n}) for $k=1$ looks as
\begin{equation}
\label{cup1}
\begin{array}{l}
dE_{1,1}(a;b)+E_{1,1}(da;b)+E_{1,1}(a;db)= a\cdot b+b\cdot a.
\end{array}
\end{equation}
So the operation $E_{1,1}$ is a sort of $\smile_1$ product: it is
the chain homotopy which which measures the noncommutativity of
$A$, c.f. the condition (\ref{Steen}). Below we denote
$a\smile_1b=E_{1,1}(a;b) $.

The condition (\ref{E2n}) for $k=1$ looks as
\begin{equation}
\label{lhirsch}
\begin{array}{l}
(a\cdot b)\smile _1c+a\cdot (b\smile _1c)+(a\smile _1c)\cdot b=0,
\end{array}
\end{equation}
this means, that the operation $E_{1,1}=\smile _1$ satisfies the
{\it left Hirsch formula} (\ref{Hirsch}).

The condition (\ref{E1n}) for $k=2$ looks as
\begin{equation}
\label{rhirsch}
\begin{array}{l}
dE_{1,2}(a;b,c)+E_{1,2}(da;b,c)+ E_{1,2}(a;db,c)+E_{1,2}(a;b,dc) \\
=a\smile _1(b\cdot c)+(a\smile _1b)\cdot c+ b\cdot (a\smile _1c),
\end{array}
\end{equation}
this means, that this $\smile _1$ satisfies the \emph{ right
Hirsch formula} just up to homotopy and the appropriate homotopy
is the operation $E_{1,2}$.

The condition (\ref{E1assoc}) for $n=m=2$ looks as
\begin{equation}
\label{cup1assoc}
\begin{array}{l}
(a\smile _1b)\smile  _1c + a\smile _1(b\smile _1c)=
E_{1,2}(a;b,c)+E_{1,2}(a;c,b),
\end{array}
\end{equation}
this means that the same operation $E_{1,2}$ measures also the
deviation from the associativity of the operation $E_{1,1}=\smile
_1$.

\subsection{hGa as a $B(\infty)$-algebra}

The notion of $B_\infty -$algebra was introduced in  \cite{GJ94}
as an additional structure on a dg module $(A,d)$ which turns the
tensor coalgebra $T^c(s^{-1}A)$ into a dg bialgebra. So it
requires a differential
$$
\widetilde{d}:T^c(s^{-1}A)\rightarrow T^c(s^{-1}A)
$$
which is a coderivation (that is an $A(\infty)$-algebra structure
on $A$, see bellow) and a an associative multiplication
$$
\widetilde{\mu }:T^c(s^{-1}A)\otimes T^c(s^{-1}A)\rightarrow
T^c(s^{-1}A)
$$
which is a map of dg coalgebras.

Here we show that a hGa structure on $A$ is a particular
$B(\infty)$-algebra structure: it induces on
$B(A)=(T^c(s^{-1}A),d_B)$ a multiplication but does not change the
differential $d_B$ (see \cite{GJ94}, \cite{Kad04}, \cite{Kad05},
\cite{KS05} for more details).

Let us extend our sequence $\{E_{1,k},\ k=0,1,2,...\}\}$ to the
sequence $\{E_{p,q}:(A^{\otimes p})\otimes (A^{\otimes q})\to A,\
p,q=0,1,...\}$ adding
\begin{equation}
\label{unit} E_{0,1}=id,\ E_{0,q>1}=0,\ E_{1,0}=id,\ E_{p>1,0}=0,
\end{equation}
and $E_{p>1,q}=0$.

This sequence defines a map $E:B(A)\otimes B(A)\to A $ by
$E([a_1,...,a_m]\otimes
[b_1,...,b_n])=E_{p,q}(a_1,...,a_m;b_1,...,b_n)$. The conditions
(\ref{E1n}) and (\ref{E2n}) mean exactly $dE+E(d_B\otimes
id+id\otimes d_B)=E\smile E$, i.e. $E$ is a twisting cochain. Thus
according to the section \ref{twistcoch} it's coextesionis a dg
coalgebra map
$$
\mu_E:B(A)\otimes B(A)\to B(A).
$$
The condition (\ref{E1assoc}) can be rewritten as $E(\mu_E\otimes
id-id\otimes \mu_E)=0, $ so this condition means that the
multiplication $\mu_E$ is associative. And the condition
(\ref{unit}) means that $[\ ]\in B(A)$ is the unit for this
multiplication.

Finally we obtained that $(B(A),d_B,\Delta,\mu_E)$ is a dg
bialgebra thus a hGa is a $B(\infty)$-algebra.


Let us mention, that a twisting cochain $E$ satisfying just the
starting condition (\ref{unit}) was constructed in \cite{Khe79}
using acyclic models for $A=C^{*}(X)$, the singular cochain
complex of a topological space. The condition (\ref{unit})
determines this twisting cochain $E$ uniquely up to equivalence of
twisting cochains (\ref{berik}).


\subsection{Homology of a hGa is a Gerstenhaber algebra}

A structure of a hGa on $A$ induces on the homology $H(A)$ a
structure of Gerstenhaber algebra (G-algebra).

Gerstenhaber algebra (see \cite{Ger63}, \cite{GV95}, \cite{Vor99})
is defined as a commutative graded algebra $(H,\cdot )$ together
with a Lie bracket of degree -1
$$
[\ ,\ ]:H^p\otimes H^q\rightarrow H^{p+q-1}
$$
(i.e. a graded Lie algebra structure on the desuspension
$s^{-1}H$) which is a biderivation: $[a,b\cdot c]=[a,b]\cdot
c+b\cdot [a,c]$. Main example of Gerstenhaber algebra is
Hochschild cohomology of an associative algebra.

The following argument shows the existence of this structure on
the homology $H(A)$ of a hGa.

First, there appears on the desuspension $s^{-1}A$  a structure of
dg Lie algebra: although the $\smile _1=E_{1,1}$ is not
associative, the condition (\ref{cup1assoc}) implies the
pre-Jacobi identity
$$
a\smile _1(b\smile _1c)+(a\smile _1b)\smile _1c= a\smile
_1(c\smile _1b)+(a\smile _1c)\smile _1b,
$$
this condition guarantees that the commutator $[a,b]=a\smile _1b
+b\smile _1a$ satisfies the Jacobi identity, besides the condition
(\ref{cup1}) implies that $[\ ,\ ]:s^{-1}A\otimes
s^{-1}A\rightarrow s^{-1}A$ is a chain map. Consequently there
appears on $s^{-1}H(A)$  a structure of graded Lie algebra. The
Hirsch formulae (\ref{lhirsch}) and (\ref{rhirsch}) imply that the
induced Lie bracket is a biderivation.

\subsection{Operadic Description}
\label{operad}

The operations $E_{1,k}$ forming hGa have nice description in
terms of the {\it surjection operad}, see \cite{MS99},
\cite{MS01}, \cite{BF01} for definition. Namely, to the dot
product corresponds the element $(1,2)\in \chi_0(2)$; to
$E_{1,1}=\smile_1$ product corresponds $(1,2,1)\in \chi_1(2)$, and
generally to the operation $E_{1,k}$ corresponds the element
\begin{equation}
\label{surj} E_{1,k}=(1,2,1,3,...,1,k,1,k+1,1)\in \chi_k(k+1).
\end{equation}
We remark  here that the defining conditions of a hGa (\ref{E10}),
(\ref{E1n}), (\ref{E2n}), (\ref{E1assoc}) can be expressed in
terms of operadic structure (differential, symmetric group action
and composition product) and the elements (\ref{surj}) satisfy
these conditions \emph{already in the operad} $\chi$.

Note that the elements (\ref{surj}) together with the element
$(1,2)$ generate the suboperad $F_2\chi$ which is equivalent to
the little square operad (\cite{MS99}, \cite{MS01}, \cite{BF01}).
This in particular implies that a cochain complex $(A,d)$ is a hGa
if and only if it is an algebra over the operad $F_2\chi$.

 This fact and the hGa structure on the
Hochschild cochain complex $C^*(U,U)$  of an algebra $U$
\cite{Kad88} were used by some authors to prove the \emph{ Deligne
conjecture} about the action of the little square operad on on the
Hochschild cochain complex $C^*(U,U)$.


\subsection{Hochschild Cochain Complex as a hGa}
\label{hoch}

Let $A$ be an algebra and $M$ be a two sided module on $A$. The
Hochschild cochain complex $C^{*}(A;M)$ is defined as
$C^n(A;M)=Hom(A^{\otimes ^n},M)$ with differential $\delta
:C^{n-1}(A;M)\to C^n(A;M)$ given by
$$
\begin{array}{ll}
\delta f(a_1\otimes ...\otimes a_n)=&a_1\cdot f(a_2\otimes ...\otimes a_n) \\
&+\sum_{k=1}^{n-1}f(a_1\otimes ...\otimes a_{k-1}\otimes a_k\cdot
a_{k+1}\otimes ..\otimes a_n)\\
&+f(a_1\otimes ...\otimes a_{n-1})\cdot a_n.
\end{array}
$$
We focus on the case $M=A$.

In this case the Hochschild complex becomes a dg algebra with
respect to the $\smile $ product defined in \cite{Ger63} by
$$
f\smile g(a_1\otimes ...\otimes a_{n+m})=f(a_1\otimes ...\otimes
a_n)\cdot g(a_{n+1}\otimes ...\otimes a_{n+m}).
$$

In \cite{Kad88} (see also \cite{GJ94}, \cite{GV95}) there are
defined the operations
$$
E_{1,i}:C^n(A;A)\otimes C^{n_1}(A;A)\otimes...\otimes
C^{n_i}(A;A)\to C^{n+n_1+...+n_i-i}(A;A)
$$
given by $E_{1,i}(f;g_1,...,g_i)=0$ for $i>n$ and
\begin{equation}
\label{cup1gg}
\begin{array}{ll}
E_{1,i}(f;g_1,...,g_i)(a_1\otimes ...\otimes
a_{n+n_1+...+n_i-i})&\\
= \sum_{k_1,...,k_i} f(a_1\otimes ...\otimes a_{k_1}\otimes
g_1(a_{k_1+1}\otimes ... \otimes a_{k_1+n_1})\otimes
a_{k_1+n_1+1}\otimes... &\\
\ \ \ \otimes a_{k_2}\otimes g_2(a_{k_2+1}\otimes ... \otimes
a_{k_2+n_2})\otimes a_{k_2+n_2+1}
\otimes ...\\
\ \ \ \otimes a_{k_i}\otimes g_i(a_{k_i+1}\otimes ... \otimes
a_{k_i+n_i})\otimes a_{k_i+n_i+1}\otimes...\otimes
a_{n+n_1+...+n_i-i}).
\end{array}
\end{equation}
The straightforward verification shows that the collection
$\{E_{1,k}\}$ satisfies the conditions  (\ref{E10}), (\ref{E1n}),
(\ref{E2n}) and (\ref {E1assoc}), thus it forms on the Hochschild
complex $C^{*}(A;A)$ a structure of homotopy G-algebra.

We note that the operation $E_{1,1}$ coincides with the circle
product defined by Gerstenhaber in \cite{Ger63}, note also that
the operation $E_{1,2}$ satisfying (\ref{rhirsch}) and
(\ref{cup1assoc}) also is defined there.


\section{Twisting Elements}
\label{twisting}

In this section we present an analog of the notion of twisting
element (see the introduction) in a homotopy G-algebra replacing
in the defining equation $da=a\cdot a$ the dot product by the
$\smile_1$ product. The appropriate notion of equivalence also
will be introduced.



Let $(C^{*,*},d,\cdot ,\{E_{1,k}\})$ be a \emph{ bigraded}
homotopy G-algebra. That is  $(C^{*,*},\cdot )$ is a bigraded
algebra $C^{m,n}\cdot C^{p,q}\subset C^{m+p,n+q} $, and we require
the existence of a differential (derivation) $d(C^{m,n})\subset
C^{m+1,n}$ and of a sequence of operations
$$
E_{1,k}:C^{m,n}\otimes C^{p_1,q_1}\otimes ...\otimes
C^{p_k,q_k}\to C^{m+p_1+...+p_k-k,n+q_1+...+q_k}
$$
so that the \emph{     total complex} (the total degree of
$C^{p,q}$ is $p+q$) is a hGa.

Bellow we introduce two versions of the notion of \emph{ twisting
elements} in a bigraded homotopy G-algebra. Although it is
possible to reduce them to each other by changing gradings, we
prefer to consider them separately in order to emphasize different
areas of their applications. The first one controls Stasheff's
$A_\infty $-deformation of graded algebras and the second controls
Gerstenhaber's deformation of associative algebras, see the next
two sections.

\subsection{Twisting Elements in a Bigraded Homotopy G-algebra (version 1)}
\label{twisthGa1}

 A \emph{ twisting element} in $C^{*,*}$ we
define as
$$
m=m^3+m^4+...+m^p+...\  ,\  m^p\in C^{p,2-p}
$$
satisfying the condition $dm=E_{1,1}(m;m)$ or changing the
notation $dm=m\smile _1m$. This condition can be rewritten in
terms of components as
\begin{equation}
\label{htwist}dm^p=\sum_{i=3}^{p-1}m^i\smile _1m^{p-i+2}.
\end{equation}
Particularly $dm^3=0,\ dm^4=m^3\smile _1m^3,\ dm^5=m^3\smile
_1m^4+m^4\smile _1m^3,...\ $. The set of all twisting elements we
denote by $Tw(C^{*,*})$.

Consider the set $G=\{g=g^2+g^3+...+g^p+...;\  g^p\in
C^{p,1-p}\}$, and let us introduce on $G$ the following operation
\begin{equation}
\label{group}\overline{g}*g=\overline{g}+g+\sum_{k=1}^\infty
E_{1,k} (\overline{g};g,...,g),
\end{equation}
particularly
$$
\begin{array}{l}
(\overline{g}*g)^2=\overline{g}+g^2;\\
(\overline{g}*g)^3=\overline{g}^3+g^3+\overline{g}^2\smile _1g^3; \\
(\overline{g}*g)^4=\overline{g}^4+g^3+
\overline{g}^2\smile _1g^3+\overline{g}^3\smile _1g^2+E_{1,2}(\overline{g}%
^2;g^2,g^2).
\end{array}
$$
It is possible to check, using the defining  conditions of a  hGa
(\ref{E10}), (\ref{E1n}), (\ref{E2n}), (\ref{E1assoc}) that this
operation is associative, has the unit $e=0+0+...$ and the
opposite $g^{-1}$ can be solved inductively from the equation
$g*g^{-1}=e$. Thus $G$ is a group.

The group $G$ acts on the set $Tw(C^{*,*})$ by the rule
$g*m=\overline{m}$ where
\begin{equation}
\label{actionm} \overline{m}=m+dg+g\cdot
g+E_{1,1}(g;m)+\sum_{k=1}^\infty E_{1,k}(\overline{m};g,...,g),
\end{equation}
particularly%
$$
\begin{array}{l}
\overline{m}^3=m^3+dg^2;\\
\overline{m}^4=m^4+dg^3+g^2\cdot g^2+g^2\smile
_1m^3+\overline{m^3}\smile _1g^2; \\
\overline{m}^5=m^5+dg^4+g^2\cdot g^3+g^3\cdot
g^2+g^2\smile _1m^4+g^3\smile _1m^3+ \\
\ \ \ \ \ \ \ \ \overline{m}^3\smile _1g^3+\overline{m}%
^4\smile _1g^2+E_{1,2}(\overline{m}^3;g^2,g^2).
\end{array}
$$
Note that although in the right hand side of the formula
(\ref{actionm}) participates $\overline{m}$ but it has less
dimension then the left hand side $\overline{m}$, thus this action
is well defined: the components of $\overline{m}$ can be solved
from this equation inductively. It is possible to check that the
resulting $\overline{m}$ is a twisting element. By $D(C^{*,*})$ we
denote the set of orbits $Tw(C^{*,*})/G.$

This group action allows us to \emph{ perturb} twisting elements
in the following sense. Let $g^n\in C^{n,1-n}$ be an arbitrary
element, then for $ g=0+...+0+g^n+0+... $ the twisting element
$\overline{m}=g*m$ looks as
\begin{equation}
\label{pert}
\overline{m}=m^3+...+m^n+(m^{n+1}+dg^n)+\overline{m}^{n+2}+\overline{m}%
^{n+3}+...\quad ,
\end{equation}
so the components $m^3,...,m^n$ remain unchanged and
$\overline{m}^{n+1}=m^{n+1}+dg^n$.

The perturbations allow to consider the following two problems.

\noindent {\bf Quantization.} Let us first mention that for a
twisting element $m=\sum m^k$ the first component $m^3\in
C^{3,-1}$ is a cycle and any perturbation does not change it's
homology class $[m^3]\in H^{3,-1}(C^{*,*})$. Thus we have the
correct map $\phi:D(C^{*,*})\to H^{3,-1}(C^{*,*})$.

A \emph{quantization} of a homology class $\alpha \in
H^{3,-1}(C^{*,*})$ we define as a twisting element $m=m^3+m^4+...$
such that $[m^3]=\alpha$. Thus $\alpha$ is quantizable if it
belongs to the image of $\phi$.

The obstructions for quatizability lay in homologies
$H^{n,3-n}(C^{*,*}),\ n\geq 5$. Indeed, let $m^3\in C^{3,-1}$ be a
cycle from $\alpha$. The first step to quantize $\alpha$ is to
construct $m^4$ such that $dm^4=m^3\smile_1m^3$. The necessary and
sufficient condition for this is $[m^3\smile_1 m^3]=0\in
H^{5,-2}(C^{*,*})$, so this homology class is the first
obstruction $O(m^3)$. Suppose it vanishes; so there exists $m^4$.
Then it is easy to see that $m^3\smile_1 m^4+m^4\smile_1 m^3$ is a
cycle and its class $O(m^3,m^4)\in H^{6,-3}(C^{*,*})$ is the
second obstruction. If $O(m^3,m^4)=0$ then there exists $m^5$ such
that $dm^5=m^3\smile_1 m^4+m^4\smile_1 m^3$. If not then we take
another $m^4$ and try new second obstruction (we remark that
changing of $m^3$ makes no sense). Generally the obstruction is
$$
O(m^3,m^4,...,m^{n-2})=[\sum_{k=3}^{n-2}m^{k}\smile_1
m^{n-k+1}]\in H^{n,3-n}(C^{*,*}).
$$

\noindent{\bf Rigidity.} A twisting element
$m=m^3+m^4+...+m^p+...$ is called \emph{trivial} if it is
equivalent to $0$. A bigraded hGa $C^{*,*}$ is \emph{rigid} if
each twisting element is trivial, i.e. if $D(C^{*,*})=\{0\}$. The
obstructions to triviality of a twisting element lay in homologies
$H^{n,2-n}(C^{*,*}),\ n\geq 3$. Indeed, for a twisting element
$m=m^3+m^4+...+m^p+...$ the first component $m^3$ is a cycle and
by (\ref{actionm}) each perturbation of $m$ leaves the class
$[m^3]\in H^{3,-1}(C^{*,*})$ unchanged and this class is the first
obstruction for triviality. If this class is zero, then we choose
$g^2\in C^{2,-1}$ such that $dg^2=m^3$. Perturbing $m$ by
$g=g^2+0+0+...$ we kill the first component $m^3$, i.e. we get the
twisting element $\overline{m}\sim m$, which looks as $
\overline{m}=0+\overline{m}^4+\overline{m}^5+...\ . $ Now, since
of (\ref{htwist}), the component $\overline{m}^4$ becomes a cycle
and it's homology class is the second obstruction. If this class
is zero then we can kill $\overline{m}^4$. If it is not then we
take another $g^2$ and try new second obstruction. Generally after
killing first components, for $m=0+0+...+0+m^{n}+m^{n+1}+...$ the
obstruction is the homology class $[m^n]\in H^{n,2-n}(C^{*,*})$.

This in particular implies that if for a bigraded homotopy
G-algebra $C^{*,*}$ all homology modules $H^{n,2-n}(C^{*,*})$ are
trivial for $n\geq 3$, then $D(C^{*,*})=0$, thus $C^{*,*}$ is
rigid.




\subsection{Twisting Elements in a Bigraded Homotopy G-algebra (version 2)}
\label{twisthGa2}

This version can be obtained from the previous one by changing
grading: take new bigraded module $\overline{C}^{p,q}=C^{p+q,-q}$.
The same operations turn $\overline{C}^{*,*}$ into a bigraded hGa.

A twisting element $m\in C^{*,2-*}$ in this case looks as
$b=b_1+b_2+...+b_n+...\ ,\  b_n\in \overline{C}^{2,n}$ where
$b_k=m^{k-2}$ and satisfies the condition $db=b\smile_1b$, or
equivalently $db_n=\sum_{i=2}^{n-1}b_i\smile _1b_{n-i}$.

Here we have the group $G^{\prime }=\{g=g_1+g_2+...+g_p+...\
;\quad g_p\in B^{1,p}\}$ with operation $ g^{\prime }*g=g^{\prime
}+g+\sum_{k=1}^\infty E_{1,k}(g^{\prime };g, ..,g)$.
This group acts on the set $Tw^{\prime }(\overline{C}^{*,*})$ by the rule $%
g*b=b^{\prime }$ where
\begin{equation}
\label{actionb} b^{\prime }=b+dg+g\cdot
g+E_{1,1}(g;b)+\sum_{k=1}^\infty E_{1,k}(b^{\prime };g,...,g).
\end{equation}
By $D^{\prime }(\overline{C}^{*,*})$ we denote the set of orbits
$Tw^{\prime }(\overline{C}^{*,*})/G^{\prime }.$


We consider the following two problems.

\noindent {\bf Quantization.} The first component $b_1\in
\overline{C}^{2,1}$ of a twisting element $b=\sum b_i$ is a cycle
and any perturbation does not change it's homology class $\alpha
=[b_1]\in H^{2,1}(\overline{C}^{*,*})$. Thus we have a correct map
$\psi:D'(\overline{C}^{*,*})\to H^{2,1}(\overline{C}^{*,*})$.

A \emph{quantization} of a homology class $\alpha \in
H^{2,1}(\overline{C}^{*,*})$ we define as a twisting element
$b=b_1+b_2+...$ such that $[b_1]=\alpha$. Thus $\alpha$ is
quantizable if $\alpha\in Im\ \psi$.

The argument similar to above shows that the obstructions to
quatizability lay in homologies $H^{3,n}(\overline{C}^{*,*}),\
n\geq 2$.

\noindent{\bf Rigidity.} A twisting element $b=b_1+b_2+...$ is
called \emph{trivial} if it is equivalent to $0$. A bigraded hGa
$\overline{C}^{*,*}$ is \emph{rigid} if each twisting element is
trivial, i.e. if $D'(\overline{C}^{*,*})=\{0\}$. The obstructions
to triviality of a twisting element lay in homologies
$H^{2,n}(\overline{C}^{*,*}),\ n\geq 1$.

This in particular implies that if for a bigraded hGa
$\overline{C}^{*,*}$ we have $H^{2,n}(\overline{C}^{*,*})=0,\
n\geq 1$ then $D'(\overline{C}^{*,*})=0$ thus $\overline{C}^{*,*}$
is rigid.



\subsection{Twisting Elements in a dg Lie Algebra corresponding to a hGa}
\label{twlie}

As it is described above for a homotopy G-algebra $(C,\cdot
,d,\{E_{1k}\})$ the desuspension $s^{-1}A$ is a dg Lie algebra
with the bracket $[a,b]=a\smile _1b-b\smile _1a$. Note that if
$C^{*,*}$ is a bigraded homotopy G-algebra, then
$L^{*,*}=s^{-1}C^{*,*}C^{*-1,*}$ is a bigraded dg Lie algebra.

Suppose $m\in C^{*,2-*}$ is a twisting element in $C^{*,*}$. The
defining equation $dm=m\smile _1m$ can be rewritten in terms of
bracket as $dm=\frac {1}{2}[m,m]$, so the same $m$ can be regarded
as a Lie twisting element in the bigraded dg Lie algebra
$L^{*,*}$.

So the notion of the twisting element in a hGa, which involves
just the operation $E_{1,1}=\smile_1$ in fact can be expressed in
terms of Lie bracket $[\ ,\ ]$.

But it is unclear whether the group action formulas
(\ref{actionm}) and (\ref{actionb}), which involve all the
operations $\{E_{1,k},\ k=1,2,...\}$ can be expressed just in
terms of bracket.


\section{ Deformation of Algebras}
\label{gerdef}

This is just illustrative application. Using the homotopy
G-algebra structure, the notions of twisting element and their
transformation one can obtain the well known results of
Gerstenhaber from \cite{Ger64}.

Let $(A,\cdot )$ be an algebra over a field $k$, $k[[t]]$ be the
algebra of formal power series in variable $t$ and
$A[[t]]=A\otimes k[[t]]$ be the algebra of formal power series
with coefficients from $A$.

Gerstenhaber's deformation of an algebra $(A,\cdot )$ is defined
as a sequence of
homomorphisms%
$$
B_i:A\otimes A\rightarrow A,\quad i=0,1,2,...;\quad B_0(a\otimes
b)=a\cdot b
$$
satisfying the \emph{ associativity} condition
\begin{equation}
\label{deform}\sum_{i+j=n}B_i(a\otimes B_j(b\otimes
c))=\sum_{i+j=n}B_i(B_j(a\otimes b)\otimes c).
\end{equation}

Such a sequence determines the \emph{ star product}
$$
a\star b=a\cdot b+B_1(a\otimes b)t+B_2(a\otimes b)t^2+B_3(a\otimes
b)t^3+...\in A[[t]],
$$
which can be naturally extended to a $k[[t]]-$bilinear product $
\star :A[[t]]\otimes A[[t]]\rightarrow A[[t]] $ and the condition
(\ref{deform}) guarantees it's associativity.

Two deformations $\{B_i\}$ and $\{B_i^{\prime }\}$ are called
\emph{ equivalent} if there exists a sequence of homomorphisms
$\{G_i:A\rightarrow A;\quad i=0,1,2,...;\quad G_0=id\}$ such that
\begin{equation}
\label{gauge}\sum_{r+s=n}G_r(B_s(a\otimes
b))=\sum_{i+j+k=n}B_i^{\prime }(G_j(a)\otimes G_k(b)).
\end{equation}
The sequence $\{G_i\}$ determines homomorphism $ G=\sum
G_it^i:A\rightarrow A[[t]]. $ On it's turn this $G$ naturally
extends to a $k[[t]]-$linear bijection $(A[[t]],\star )\rightarrow
(A[[t]],\star ^{\prime })\,$ and the condition (\ref{gauge})
guarantees that this extension is multiplicative.

A deformation $\{B_i\}$ is called \emph{ trivial}, if $\{B_i\}$ is
equivalent to $\{B_0,0,0,...\}$.  An algebra $A$ is called \emph{
rigid}, if each it's deformation is trivial.

Now we present the interpretation of deformations and their
equivalence in terms of twisting elements of version 2 type and
their equivalence in hGa of Hochschild cochains.

As it is mentioned in \ref{hoch} the Hochschild complex
$C^{*}(A,A)$ for an algebra $A$ is a homotopy G-algebra. Then the
tensor product $ C^{*,*}=C^{*}(A,A)\otimes k[[t]] $ is a \emph{
bigraded }Hirsch algebra with the following structure
$$
\begin{array}{c}
C^{p,q}=C^p(A,A)\cdot t^q,\ \delta(f\cdot t^q)=\delta f\cdot t^q,\
\
f\cdot t^p\smile g\cdot t^q=(f\smile g)\cdot t^{p+q},\\
E_{1,k}(f\cdot t^p;g_1\cdot t^{q_1}, ..., g_k\cdot
t^{q_k})=E_{1,k}(f;g_1, ..., g_k)\cdot t^{p+q_1+...+q_k},
\end{array}
$$
here we use the notation $f\otimes t^p=f\cdot t^p$.

Then each deformation $\{B_i:A^{\otimes ^2}\to A,\ i=1,2,3,...\}$
can be interpreted as a version 2 type twisting element
$b=b_1+b_2+...+b_k+...,\ b_k=B_k\cdot t^k\in C^{2,k}$: the
associativity condition (\ref{deform}) can be rewritten as%
$$
\delta B_n\cdot t^n=\sum_{i+j=n}B_i\cdot t^i\smile _1B_j\cdot t^j.
$$

Suppose now two deformations $\{B_i\}$ and $\{B'_i\}$ are
equivalent, i.e. there exists $\{G_i\}$ such that the condition
(\ref{gauge}) satisfied. In terms of Hochschild cochsins this
condition looks as
$$
b'=b+\delta g+g\smile g+g\smile _1b+E_{1,1}(b';g)+E_{1,2}(b';g,g),
$$
where $g=g_1+...+g_k+...,\ g_k=G_k\cdot t^k\in C^{1,k}$. This
equality slightly differs from (\ref{actionb}), but since $g\in
C^1(A,A)$ and $b'\in C^2(A,A))$, we have $E_{1,k}(b';g,...,g)=0$
for $k\geq 3$ (see \ref{hoch}), thus they in fact coincide.

 So we obtain that
deformations are equivalent if and only if the corresponding
Hochschild twisting elements $b$ and $b'$ are equivalent.
Consequently the set of equivalence classes of deformations is
bijective to $D'(C^{*,*})$.

It is clear that $H^{p,q}(C^{*,*})=HH^p(A,A)\cdot t^q$. Then from
the section \ref{twisthGa2} follow the classical results of
Gerstenhaber: obstructions for quantization of a homomorphism
$b_1:A\otimes A\to A$ lay in $HH^3(A,A)$, and if $HH^3(A,A)=0$
then each $b_1$ is quantizable (or \emph{integrable} as it is
called in \cite{Ger64}). Furthermore, the obstructions for
triviality of a deformation lay in $HH^2(A,A)$, and if
$HH^2(A,A)=0$ then $A$ is rigid.

\begin{remark}
As we see in the definition of equivalence of deformations
participate just the operations $E_{1,1}$ and $E_{1,2}$, the
higher operations $E_{1,k},\ k>2$ disappear since of
(\ref{cup1gg}). So observing just deformation problem it is
impossible to establish general formula (\ref{actionb}) for
transformation of twisting elements.
\end{remark}


\section{$A(\infty )$-deformation of Graded Algebras}
\label{stadef}

In this section we give the similar description of $A(\infty
)$-deformation of graded algebras in terms of twisting elements in
the hGa of Hochschild cochains. So this two types of deformation
will be unified by the notion of twisting element in hGa.
Partially these results are given in \cite{Kad88}, \cite{Kad93}.

\subsection{$A(\infty )$-algebras}

The notion of $A(\infty )$-algebra was introduced by J.D. Stasheff
in \cite {Sta63}. This notion generalizes the notion of dg
algebra.

An $A(\infty )$-algebra is a graded module $M$ with a given
sequence of
operations%
$$
\{m_i:M^{\otimes i}\rightarrow M,\quad i=1,2,...,\quad \deg
m_i=2-i\}
$$
which satisfies the following conditions
\begin{equation}
\label{Ainf} \sum_{i+j=n+1}\sum_{k=0}^{n-j}m_i(a_1\otimes
...\otimes a_k\otimes m_j(a_{k+1}\otimes ...\otimes
a_{k+j})\otimes ...\otimes a_n)=0.
\end{equation}
Particularly, for the operation $m_1:M\rightarrow M$ we have $\deg
m_1=1$ and $m_1m_1=0$, this $m_1$ can be regarded as a
differential on $M$. The
operation $m_2:M\otimes M\rightarrow M$ is of degree 0 and satisfies%
$$
m_1m_2(a_1\otimes a_2)+m_2(m_1a_1\otimes a_2)+m_2(a_1\otimes
m_1a_2)=0,
$$
i.e. $m_2$ can be regarded as a multiplication on $M$ and $m_1$ is
a derivation. Thus $(M,m_1,m_2)$ is a sort of (maybe
nonassociative) dg algebra. For the operation $m_3$ we have $\deg
m_3=-1$ and
$$
\begin{array}{l}
m_1m_3(a_1\otimes a_2\otimes a_3)+m_3(m_1a_1\otimes a_2\otimes
a_3)+m_3(a_1\otimes m_1a_2\otimes a_3)\\
+ m_3(a_1\otimes a_2\otimes m_1a_3)+m_2(m_2(a_1\otimes a_2)\otimes
a_3)+m_2(a_1\otimes m_2(a_2\otimes a_3))=0,
\end{array}
$$
thus the multiplication $m_2$ is \emph{ homotopy associative} and
the appropriate chain homotopy is $m_3$.

The sequence of operations $\{m_i\}$ determines on the tensor
coalgebra
$$
T^c(s^{-1}M)=R +s^{-1}M+s^{-1}M\otimes s^{-1}M+s^{-1}M\otimes
s^{-1}M\otimes s^{-1}M+...
$$
a coderivation
$$
d_m(a_1\otimes ...\otimes a_n)=\sum_{k,j}a_1\otimes ...\otimes
a_k\otimes m_j(a_{k+1}\otimes ...\otimes a_{k+j})\otimes
...\otimes a_n,
$$
and the condition (\ref{Ainf}) is equivalent to $d_md_m=0$. The
obtained dg coalgebra $(T^c(s^{-1}M),d_m)$ is called \emph{ bar
construction } and is denoted as $B(M,\{m_i\})$.

A morphism of $A(\infty )$-algebras $(M,\{m_i\})\rightarrow
(M^{\prime },\{m_i^{\prime }\})$ is defined as a sequence of
homomorphisms
$$
\{f_i:M^{\otimes ^i}\rightarrow M^{\prime },\quad i=1,2,...,\quad
\deg f_i=1-i\},
$$
which satisfy the following condition
\begin{equation}
\label{Ainfmor}
\begin{array}{l}
\sum_{i+j=n+1}\sum_{k=0}^{n-j}f_i(a_1\otimes ...\otimes a_k\otimes
m_j(a_{k+1}\otimes ...\otimes a_{k+j})\otimes ...\otimes a_n)\\
= \sum_{k_1+...+k_t=n} m'_t(f_{k_1}(a_1\otimes...\otimes
a_{k_1})\otimes f_{k_2}(a_{k_1+1}\otimes... \otimes
a_{k_1+k_2})\\
\ \ \ \otimes ...\otimes
f_{k_t}(a_{k_1+...+k_{t-1}+1}\otimes...\otimes a_{n})).
\end{array}
\end{equation}
In particular for $n=1$ this condition gives
$f_1m_1(a)=m'_1f_1(a)$, i.e. $f_1:(M,m_1)\to (M',m'_1)$ is a chain
map; for $n=2$ it gives
$$
\begin{array}{l}
f_1m_2(a_1\otimes a_2)+m'_2(f_1(a_1)\otimes f_1(a_2))\\
=m'_1f_2(a_1\otimes a_2)+ f_2(m_1a_1\otimes a_2)+f_2(a_1\otimes
m_1a_2),
\end{array}
$$
thus $f_1:(M,m_1,m_2)\to (M',m'_1,m'_2)$ is multiplicative up to
homotopy $f_2$.

A  collection $\{f_i\}$ defines a homomorphism $f:B(M,\{m_1\})\to
M'$. It's \emph{ comultiplicative coextension}, see \ref{bar}, is
a graded coalgebra map of the bar constructions
$$
B(f):B(M,\{m_i\})\rightarrow B(M^{\prime },\{m_i^{\prime }\}),
$$
and the condition (\ref{Ainfmor}) guarantees that $B(f)$ is a
cahin map, i.e. $B(f)$ is a morphism of dg coalgebras. So $B$ is a
functor from the category of $A(\infty)$-algebras to the category
of dg coalgebras.

A weak equivalence of $A(\infty) $-algebras is defined as a
morphism $f=\{f_i\}$ for which $B(f)$ is a weak equivalence of dg
coalgebras. It is possible to show (see for example \cite{Kad93})
that

\noindent (i) $f$ is a weak equivalence of $A(\infty )$-algebras
if and only if $f_1$ is a weak equivalence of dg modules;

\noindent (ii) $f$ is an isomorphism of $A(\infty )$-algebras if
and only if $f_1$ is an isomorphism of dg modules.

An $A(\infty )$-algebra $(M,\{m_i\})$ we call \emph{ minimal} if
$m_1=0$. In this case $(M,m_2)$ is \emph{ strictly} associative
graded algebra.

The following proposition is the imediate consequence of (i) and
(ii):
\begin{proposition}
Each weak equivalence of minimal $A(\infty )$-algebras is an
isomorphism.
\end{proposition}

\subsection{$A(\infty )$ Deformation of Graded Algebras as
Twisting Element}

Let $(H,\mu :H\otimes H\rightarrow H)$ be a graded algebra. It's
Stasheff's (or $A(\infty )$) deformation we define as a minimal
$A(\infty)$-algebra $(H,\{m_i\})$ with $m_2=\mu$. Two deformations
$(H,\{m_i\})$ and $(H,\{m'_i\})$ we call equivalent if there
exists an isomorphism of $A(\infty)$-algebras
$\{f_i\}:(H,\{m_i\})\to (H,\{m'_i\})$ with $f_1=id$.


A deformation $(H,\{m_i\})$ we call trivial if it is equivalent to
$(H,\{m_1=0,m_2=\mu,m_{\geq 3}=0\})$. An algebra $(H,\mu)$ we call
\emph{ rigid} (or \emph{intrinsically formal}, this term is
borrowed from the rational homotopy theory), if each it's
deformation is trivial.

Now we present the interpretation of deformations and their
equivalence in terms of twisting elements and their equivalence in
hGa of Hochschild cochains.

The Hochschild cochain complex of a graded algebra $H$ with
coefficients in itself is bigraded: $C^{m,n}(H,H)=Hom^n(H^{\otimes
m},H)$, here $Hom^n$ denotes degree $n$ homomorphisms.
The coboundary operator $\delta $ maps $C^{m,n}(H,H)$ to $%
C^{m+1,n}(H,H)$. Besides, for $f\in C^{m,n}(H,H)$ and $g_i\in
C^{p_i,q_i}(H,H)$ one has $f\smile g\in C^{m+p,n+q}(H,H),\ f\smile
_1g\in C^{m+p-1,n+q}(H,H)$ and
$$
E_{1,k}(f;g_1, ..., g_k)\in
C^{m+p_1+...+p_k-k,n+q_1+...+q_k}(H,H),
$$
thus the Hochschild complex $C^{*,*}(H,H)$ is a \emph{ bigraded}
homotopy G-algebra in this case. Let us denote the $n$-the
homology module of the complex $(C^{*,k}(H,H),\delta )$ by
$HH^{n,k}(H,H)$.

Suppose now that $(H,\{m_i\})$ is an $A(\infty )$ deformation of
$H$. Each operation $m_i:H^{\otimes ^i}\rightarrow H$ can be
regarded as a Hochschild cochain from $C^{i,2-i}(H,H)$. The
condition (\ref{Ainf}) can be rewritten as
$$
\delta m_k=\sum_{i=3}^{k-1}m_i\smile _1m_{k-i+2},
$$
thus $m=m_3+m_4+...$ is a \emph{ \ twisting element} (version 1)
in $C^{*,*}(H,H)$.

Now let $(H,\{m_i\})$ and $(H,\{m_i^{\prime }\})$ be two $A(\infty
)$ deformations of $H$. Then it follows from (\ref{actionm}) that
the corresponding twisting elements $m$ and $m^{\prime }$ are
equivalent if and only if these two $A(\infty )$ deformations are
equivalent: if $m^{\prime }=p*m,$ then
$\{p_i\}:(H,\{m_i\})\rightarrow
(H,\{m_i^{\prime }\})$ with $p_1=id$ is an \emph{ isomorphism} of $A(\infty )$%
-algebras. So we obtain the
\begin{theorem}
The set of isomorphism classes of all $A(\infty)$ deformations of
a graded algebra $(H,\mu)$ is bijective to the set of equivalence
classes of twisting elements $D(C^{*,*}(H,\mu))$.
\end{theorem}
Moreover, from \ref{twisthGa1} we get the following
\begin{theorem}
\label{Ho3}If for a graded algebra $(H,\mu )$ it's Hochschild
cohomology modules $HH^{n,2-n}(H,H)$ are trivial for $n\geq 3$,
then $(H,\mu)$ is intrinsically formal.
\end{theorem}

\subsubsection{$A(\infty )$-algebra Structure in Homology of a dg algebra}

Let $(A,d,\mu )$ be a dg algebra and $(H(A),\mu ^{*})$ be it's
homology algebra. Although the product in $H(A)$ is associative,
there appears a structure of a (generally nondegenerate) minimal
$A(\infty )$-algebra, which is an $A(\infty)$ deformation of
$(H(A),\mu^*)$. Namely, in \cite{Kad76}, \cite {Kad80} the
following result was proved (see also \cite{Smi80}, \cite{GS86}):

\begin{theorem}
If for a dg algebra all homology $R$-modules $H_i(A)$ are free,
then there exist: a structure of minimal $A(\infty )$-algebra
$(H(A),\{m_i\})$ on $H(A)$ and a weak equivalence of $A(\infty
)$-algebras
$$
\{f_i\}:(H(A),\{m_i\})\rightarrow (A,\{d,\mu ,0,0,...\})
$$
such, that $m_1=0$, $m_2=\mu ^{*}$, $f_1^{*}=id_{H(A)}$, such a
structure is unique up to isomorphism in the category of $A(\infty
)$-algebras.
\end{theorem}

Particularly an $A(\infty )$-algebra structure appears in
cohomology $H^*(X)$ of a topological space $X$ or in homology
$H_*(G)$ of a topological group or H-space $G$. (Co)homology
algebra equipped with this additional structure carries more
information then just the (co)homology algebra. Some applications
of this structure are given in \cite{Kad80} , \cite{Kad93}. For
example the cohomology $A(\infty)$-algebra $(H^*(X),\{m_i\})$
determines cohomology of the loop space $H^*(\Omega X)$ when just
the algebra $(H^*(X),m_2)$ does not. Similarly, the homology
$A(\infty)$-algebra $(H_*(G),\{m_i\})$ determines homology of the
classifying space $H_*(B_G)$ when just the Pontriagin algebra
$(H_*(G),m_2)$ does not. Furthermore, the rational cohomology
$A(\infty)$-algebra $(H^*(X,Q),\{m_i\})$ (which actually is
$C(\infty)$ in this case) determines the rational homotopy type of
1-connected $X$ when just the cohomology algebra $(H^*(X,Q),m_2)$
does not.

Therefore it is of particular interest the cases, when this
additional structure is vanishes, that is when $A(\infty
)$-algebra $(H(A),\{m_i\})$ is degenerate (in this case a dg
algebra $A$ is called \emph{ formal}). The above theorem \ref{Ho3}
gives the sufficient condition of formality of $A$ in terms of
Hochschild cohomology of $H(A)$.

A. Razmadze Mathematical Institute,

1, M. Alexidze Str., Tbilisi, 0193, Georgia

kade@rmi.acnet.ge

\end{document}